\documentclass{amsart}
\title{Instantons on $S^{4}$ and $\cpbar $, rank stabilization, and Bott
periodicity. }
\author{Jim Bryan and Marc Sanders}
\date{\today}

\usepackage{amsmath,amsthm,amsfonts,amscd}

\newtheorem{thm}{Theorem}[section]
\newtheorem{cor}[thm]{Corollary}
\newtheorem{lem}[thm]{Lemma}
\newtheorem{prop}[thm]{Proposition}

\theoremstyle{definition}

\theoremstyle{remark}
\newtheorem{rem}{Remark}[section]

\newcommand{\cnums} {{\mathbf C}}          
\newcommand{\znums} {{\mathbf Z}}		
\newcommand{\cpbar} {\overline{\mathbf {CP}}^2}
\newcommand{\dbar}{\bar\partial}

\newcommand{\dirac}{\partial
\hspace{-5pt}\raisebox{1.5pt}{/}\hspace{-1pt}}

\newcommand{\Iso}{\operatorname{Iso}}
\newcommand{\Ker}{\operatorname{Ker}}

\newcommand{\Coker}{\operatorname{Coker}}
\newcommand{\im}{\operatorname{Im}}

\newcommand{\comp}{\hspace{-3pt}\stackrel{\scriptstyle
\circ}{}\hspace{-3pt}}                                 
\newcommand{\ie}{{\em i.e. }}

\newcommand{\til}[1]{\widetilde{#1}}

\newcommand{\1}{{{\mathchoice {\rm 1\mskip-4mu l} {\rm 1\mskip-4mu l}
{\rm 1\mskip-4.5mu l} {\rm 1\mskip-5mu l}}}}
\newcommand{\M}[2]{\mathcal{M}_{#1}^{#2}}
\newcommand{\E}{\mathcal{E}}
\renewcommand{\O}{\mathcal{O}}
\newcommand{\U}{A}
\renewcommand{\P}{\mathbf{P}}
\newcommand{\embed}{\hookrightarrow}

\begin{document}

\markboth{Rank Stable Instantons}{Rank Stable Instantons}
\renewcommand{\sectionmark}[1]{}
\begin{abstract}
We study the large $n$ limit of the moduli spaces of $G_{n}$-instantons on
$S^{4}$ and $\cpbar $ where $G_{n}$ is $SU(n)$, $Sp(n/2)$, or $SO(n)$. We
show that in the direct limit topology, the moduli space is homotopic to a
classifying space. For example, the moduli space of $Sp(\infty )$ or
$SO(\infty )$ instantons on $\cpbar $ has the homotopy type of $BU(k)$
where $k$ is the charge of the instantons. We use our results along with
Taubes' result concerning the $k\to \infty $ limit to obtain a novel proof
of the homotopy equivalences in the eight-fold Bott periodicity spectrum.
We work with the algebro-geometric realization of the instanton spaces as
moduli spaces of framed holomorphic bundles on $\cnums \P ^{2}$ and
$\cnums \P ^{2}$ blown-up at a point. We give explicit constructions for
these moduli spaces (see Table 1).
\end{abstract}

\maketitle

\section{Introduction}

Let $\M{k}{G_{n}}(X)$ denote the space of (based) $G_{n}$-instantons on $X$
where $G_{n}$ is $SU(n)$, $SO(n)$, or $Sp(n/2)$. In 1989, Taubes
\cite{Tau-stable} showed that there is a ``gluing'' map
$\M{k}{G_{n}}(X)\embed \M{k'}{G_{n}}(X)$ when $k'>k$.
He proved that in the direct limit topology,
the instantons
capture all the topology of connections modulo gauge equivalence. In other
words, there is a homotopy equivalence:
$$
\operatorname{lim}_{k\to \infty }\M{k}{G_{n}}(X) \sim
\operatorname{Map}_{0} (X,BG_{n}).
$$

There is also an inclusion
$\M{k}{G_{n}}(X)\embed  \M{k}{G_{n'}}(X)$ where $n'>n$ induced by the
inclusion $G_{n}\embed G_{n'}$. Not much is known
about the homotopy type of $\M{k}{G}(X)=\operatorname{lim}_{n\to \infty }\M
{k}{G_{n}}(X)$ for general $X$. In this paper we determine the homotopy
type of $\M{k}{G}(X)$ when $X$ is $S^{4}$ or $\cpbar $ with their standard
metrics. The results are

\begin{eqnarray}\label{eqn:rank stable for S4}
\M{k}{G}(S^{4})&\sim &\begin{cases}
BU(k)&\text{if }G=SU,\\
BO(k)&\text{if }G=Sp,\\
BSp(k/2)&\text{if }G=SO;
\end{cases}\\
\M{k}{G}(\cpbar )&\sim &\begin{cases}
BU(k)\times BU(k)&\text{if }G=SU,\\
BU(k)&\text{if $G$ is $Sp$ or $SO$.}
\end{cases}\nonumber
\end{eqnarray}

The results for the $S^{4}$ case were first proved in \cite{Sanders},
\cite{No-Sa}, and
\cite{Kir} (c.f. \cite{Tian}) and we proved the $\M{k}{SU}(\cpbar )$ result
in \cite{Br-Sa}. In this paper we are able to provide a unified approach to
these moduli spaces and stabilization results.

By  employing Taubes'
theorem and by utilizing the conformal map $f:\cpbar \to S^{4}$ to compare
$\M{k}{G_{n}} (S^{4 })$ and $\M{k}{G_{n}} (\cpbar )$, and we are able to
give a novel proof of
the homotopy equivalences in the real and unitary Bott periodicity
spectrums.  Work in this
direction has been  done by Tian using instantons on $S^{4}$ (see
\cite{Sanders},\cite{Tian}) where one can prove some of the 4-fold
equivalences. By using the comparison with instantons on $\cpbar $ we are
able to recover the finer 2-fold equivalences in the periodicity spectrum.

The moduli spaces $\M{k}{G_{n}}(S^{4})$ and $\M{k}{G_{n}}(\cpbar )$ are
known to be isomorphic to moduli spaces of certain holomorphic bundles and
have been constructed in various guises (\cite{BU86}, \cite{King},
\cite{DonMonads} ). Using work of Donaldson and King we construct the
spaces from a unified viewpoint (see Table 1). We
describe the relevant
moduli spaces of holomorphic bundles as follows:

Let $H\subset \cnums \P ^{2}$ be a fixed hyperplane and let
$\til{\mathbf{CP}}^{2}$  be the blow-up of $\cnums \P ^{2}$ at a
point away from $H$. Donaldson showed \cite{DonMonads} that
$\M{k}{SU(n)}(S^{4}) $ is isomorphic to the moduli space of pairs $(\E
,\tau )$ where $\E \to \cnums \P ^{2}$ is a rank $n$ holomorphic bundle with
$c_{1}(\E )=0$, $c_{2}(\E )=k$ and $\tau :\E |_{H}\to \cnums ^{n}\otimes \O
_{H}$ is a trivialization of $\E $ on $H$. In  \cite{King}, King extended
this result to
$\cpbar $ by showing that $\M{k}{SU(n)}(\cpbar ) $ is
isomorphic to the moduli space of pairs $(\E ,\tau )$ where $\E \to
\til{\cnums \P }^{2}$ is a rank $n$ holomorphic bundle with
$c_{1}(\E )=0$, $c_{2}(\E )=k$ and $\tau :\E |_{H}\to \cnums ^{n}\otimes \O
_{H}$ is a trivialization of $\E $ on $H$. They also construct the moduli
spaces in terms of ``linear algebra data.''

One can extend their results to $Sp(n/2) $ and $SO(n)$. Let $X$ denote
$S^{4}$ or $\cpbar $. The moduli space of $Sp$-instantons
(respectively $SO$-instantons)
is isomorphic to the moduli space of triples $(\E ,\tau ,\phi )$ where $\phi $
is a symplectic (resp. real) structure:

\begin{eqnarray*}
\M{k}{Sp(n/2)}(X )&\cong &\{(\E ,\tau ,\phi ):(\E,\tau )\in \M
{k}{SU(n)} (X ), \phi
:\E\stackrel{\cong }{\longrightarrow }\E ^{*}, \phi ^{*}= -\phi \},\\
\M{k}{SO(n)}(X )&\cong &\{(\E ,\tau ,\phi ):(\E,\tau ) \in \M
{k}{SU(n)} (X ),\phi
:\E\stackrel{\cong }{\longrightarrow }\E ^{*}, \phi ^{*}= \phi \}.
\end{eqnarray*}

Our construction realizes these moduli spaces as quotients of affine
varieties $A_{k}^{G_{n}}(X)$ (the ``linear algebra data'') by free actions.
The key to proving our stability theorem is to show that in the large $n$
limit, the space of ``linear algebra data'' becomes contractible.

The constructions also allow us to identify the universal bundles over
$\M{k}{G}(X)$ inducing the homotopy equivalences of Equation \ref{eqn:rank
stable for S4}. In the holomorphic setting they can be described as
certain higher direct image bundles and in the connection setting they can
be described as the index bundles of a certain family of coupled Dirac
operators.

In section \ref{sec: main result and BP} we fix notation, state the
theorems, and  prove Bott periodicity. In the subsequent sections
we construct the moduli spaces and prove the theorems. We conclude with
a short appendix discussing a more differentio-geometric construction of
the universal bundles. The authors would like to thank John Jones  and
Ralph Cohen for
suggesting that the homotopy equivalences of Equation \ref{eqn:rank stable
for S4}  should exist.

\section{The main results and Bott periodicity}\label{sec: main result and
BP}
\subsection{Statement of the theorems}\label{subsec:statement of thms}

Let $G_{n}\embed P\to X$ be a principal bundle on a Riemannian $4$-manifold
$X$ with structure group $G_{n}=SU(n)$, $SO(n)$, or $Sp(n/2)$. Using the
defining representations for $SU(n)$ or $Sp(n/2)$ and the complexified
standard representation for $SO(n)$, we associate to $P$ a rank $n$ complex
vector bundle $E$ and we define the {\em charge} $k$ to be $c_{2}(E)[X]$.
\footnote{Our definition of $k$ in the $SO(n)$ case differs from some
of the literature by a factor of $2$.}

A bundle isomorphism $\phi :E\to E^{*}$ is called a {\em real structure} if
$\phi^{*}=\phi $ and it is called a {\em symplectic structure} if $\phi
^{*} =-\phi $. We can regard a $SO(n) $ or a $Sp(n/2)$ bundle as a
$SU(n)$ bundle $E$  along with $\phi $, a real or symplectic structure
respectively.
Obviously, $n$ must be even for $E$ to have a symplectic structure and it
is also not hard to see that if $E$ has a real structure, then our $k$ must be
even.

Let $\mathcal{A}(E)$ denote the space of connections on $E$ that are
compatible with $\phi $ and let
$F^{+}_{E} $ be the self-dual part of the curvature of a connection $A\in
\mathcal{A}(E)$. Let $\mathcal{G}_{E}^{0}$ be the group of gauge
transformations of $E$ commuting with $\phi $ and preserving a fixed
isomorphism $E_{x_{0}}\cong
\cnums ^{n}$ of the fiber over a base point $x_{0}\in X$.
We define the (based) instanton moduli spaces to be (c.f. \cite{D-K}):
$$
\M{k}{G_{n}}(X)=\{A\in \mathcal{A}(E):F_{A}^{+}=0 \}/\mathcal{G}_{E}^{0}.
$$

>From here on let $X$ denote $S^{4}$ or $\cpbar $ with their standard metrics.
We describe how the moduli spaces $\M{k}{G_{n}}(X)$ can be
constructed from configurations of linear algebra data satisfying certain
``integrability'' conditions, modulo natural automorphisms. The
configurations are laid out by Table 1 where we have adopted the
following notations: Our vector spaces are always complex and
our maps are always complex linear. We regard a map $f:U\to W$ as an
element $f\in U^{*}\otimes W$. An isomorphism  $\phi :W\to  W^{*}$ is
called a {\em symplectic structure} if $\phi \in \Lambda ^{2}W^{*}$ and a
{\em real structure } if $\phi \in S^{2}W^{*}$. $Gl(W)$ denotes the group
of isomorphism of $W$ and if $\phi $ is a symplectic (respectively real)
structure on $W$, the let $Sp(W)$ (resp. $O(W)$) denote the group of
isomorphisms of $W$ compatible with $\phi $ (\ie $f^{*}\phi f=\phi $).
When $n$ is even, let $J$ denote the standard symplectic structure on
$\cnums ^{n}$. Unless otherwise noted, the vector spaces in Table
1 are $k$-dimensional.

If $f:V_{1}\to V_{2}$, then $Gl(V_{1})\times Gl(V_{2})$ acts on $f$ by
$f\mapsto g_{2}fg_{1}^{-1} $ and thus on $f^{*}$ by
$(g_{1}^{-1})^{*}f^{*}g_{2}^{*}$. So the action of the
automorphism group on $\cpbar $ configurations is given by
\footnote{Note that if we
fix bases for the vector spaces, $f^{*}$ is the transpose matrix  and
should not be confused with conjugate transpose.}
\begin{eqnarray*}
(g,h)\cdot
(a_{1},a_{2},x,b,c)&=&(ga_{1}h^{-1},ga_{2}h^{-1},hxg^{-1},gb,ch^{-1})\\
g\cdot (\alpha _{1},\alpha _{2},\xi ,\gamma )&=&(g\alpha _{1}g^{*},g\alpha
_{2}g^{*},(g^{*})^{-1}\xi g^{-1},\gamma g^{*}),
\end{eqnarray*}
and on $S^{4}$ configurations by
\begin{eqnarray*}
g\cdot (a_{1},a_{2},b,c)&=&(ga_{1}g^{-1},ga_{2}g^{-1},gb,cg^{-1})\\
g\cdot (\alpha _{1},\alpha _{2},\gamma )&=&(g\alpha _{1}g^{-1},g\alpha
_{2}g^{-1},\gamma g^{-1}).
\end{eqnarray*}

The three main theorems of this paper are the following:

\begin{table}\label{table}
\begin{picture}(500,350)
\put(-20,170){$
\begin{array}{||c|c|c|c|c|c||}\hline
G_{n}&	X&	\text{Configurations}&	\text{Integrability}&
\text{Automorphism}&\dim _{\cnums }\M{}{}\\
&&&	\text{Condition}&\text{Group}&\\ \hline \hline
SU(n)&	\cpbar &(a_{1},a_{2},x,b,c)
&a_{1}xa_{2}-a_{2}xa_{1}+bc=0&
Gl(W)\times Gl(U)&4nk\\ \cline{3-3}\cline{3-3}
&&a_{i}\in U^{*}\otimes W&&&\\
&&x\in W^{*}\otimes U&&&\\
&&b\in \cnums^{n}\otimes W&&&\\
&&c\in U^{*}\otimes \cnums^{n}&&&\\ \hline
Sp(n/2)&\cpbar &(\alpha _{1},\alpha _{2},\xi,\gamma )&
\alpha _{1}\xi \alpha _{2}-\alpha _{2}\xi \alpha _{1}-\gamma ^{*}J\gamma =0
& Gl(W)&(n+2)k\\ \cline{3-3}\cline{3-3}
&& \alpha _{i}\in S^{2}W&&&\\
&& \xi \in S^{2}W^{*}&&&\\
&& \gamma\in W\otimes \cnums ^{n}&&&\\ \hline
SO(n)&\cpbar &(\alpha _{1},\alpha _{2},\xi,\gamma )&
\alpha _{1}\xi \alpha _{2}-\alpha _{2}\xi \alpha _{1}+\gamma ^{*}\gamma =0
& Gl(W)&(n-2)k\\ \cline{3-3}\cline{3-3}
&& \alpha _{i}\in \Lambda ^{2}W&&&\\
&& \xi \in \Lambda ^{2}W^{*}&&&\\
&& \gamma\in W\otimes \cnums ^{n}&&&\\ \hline
SU(n)&	S^{4}&(a_{1},a_{2},b,c)&[a_{1},a_{2}]+bc=0&Gl(W)&4nk\\ \cline{3-3}
\cline{3-3}
&&a_{i}\in W^{*}\otimes W&&&\\
&&b\in \cnums^{n}\otimes W&&&\\
&&c\in W^{*}\otimes \cnums^{n}&&&\\ \hline
Sp(n/2)&S^{4}&(\alpha _{1},\alpha _{2},\gamma )&\Phi [\alpha _{1},\alpha
_{2}]-\gamma ^{*}J\gamma  =0&O(W)&(n+2)k\\ \cline{3-3}\cline{3-3}
&&\text{Real str. }\Phi &&&\\
&&\Phi \alpha _{i}\in S^{2}W^{*}&&&\\
&& \gamma\in W^{*}\otimes \cnums ^{n}&&&\\ \hline
SO(n)&S^{4}&(\alpha _{1},\alpha _{2},\gamma )&\Phi [\alpha _{1},\alpha
_{2}]+\gamma ^{*}\gamma  =0&Sp(W)&(n-2)k\\ \cline{3-3}\cline{3-3}
&&\text{Sympl. str. }\Phi &&&\\
&&\Phi \alpha _{i}\in \Lambda ^{2}W^{*}&&&\\
&& \gamma\in W^{*}\otimes \cnums ^{n}&&&\\ \hline
\end{array}
$}
\end{picture}
\caption{Configurations that construct $\M{k}{G_{n}}(X)$.}
\end{table}

\begin{thm}[Moduli Construction]\label{thm:moduli construction}
Let $\overline{A}^{G_{n}}_{k}(X)$ denote the space of integrable
configurations as given by Table 1. There is an open dense set
$A^{G_{n}}_{k}(X)\subset
\overline{A}^{G_{n}}_{k}(X) $ (the ``non-degenerate'' configurations) such
that the instanton moduli space $\M{k}{G_{n}}(X)$ is
isomorphic to the quotient of $A^{G_{n}}_{k}(X)$ by the automorphism group.
Furthermore, the action of the automorphism group on $A^{G_{n}}_{k}(X)$ is
free and  the vector spaces $W$ and $U$ of Table 1 are
canonically isomorphic to $H^{1}(\E (-H))$ and $H^{1}(\E (-H+E))$
respectively.
\end{thm}

\begin{thm}[Lifts of instanton maps to configurations]
\label{thm:lifts of maps to configurations}
There are commuting inclusions of configurations (defined in section
\ref{sec:lifts of maps})
\begin{eqnarray*}
i&:&A^{G_{n}}_{k}(X)\embed A^{G_{n'}}_{k}(X)\\
j&:&A^{G_{n}}_{k}(S^{4})\embed A^{G_{n}}_{k}(\cpbar )
\end{eqnarray*}
for $n<n'$ and $k<k'$. These maps intertwine the automorphisms and
consequently descend to maps on the instanton moduli spaces.  The map $i$
descends to the map induced by the inclusion $G_{n}\embed G_{n'}$ and the
map $j$ descends to the map induced by pulling back connections via
$f:\cpbar \to S^{4}$.
\end{thm}
\begin{thm}[Rank Stabilization]\label{thm:rank stabilization}
Let $A^{G}_{k}(X)$ be the direct limit space
$$\operatorname{lim}_{n\to \infty }A^{G_{n}}_{k}(X)$$
defined by the inclusions $i$. Then
$A^{G}_{k}(X)$ is contractible and consequently
$$\M{k}{G}=\operatorname{lim}_{n\to \infty }\M{k}{G_{n}}$$
is homotopic to
the classifying space for the associated automorphism group. This theorem
implies the homotopy equivalences in Equation \ref{eqn:rank stable for S4}.
\end{thm}

\begin{rem}\label{rem:dimension count}
A na\"{\i}ve dimension count for $\M{k}{G_{n}}(X)$ is obtained by subtracting
the dimension of the automorphism group and the number of conditions
imposed by integrability from the dimension of the configurations. This
agrees with the dimension predicted by the Atiyah-Singer index formula and
appears in the far right column of the table.
\end{rem}

We prove Theorems \ref{thm:moduli construction}, \ref{thm:lifts of maps to
configurations}, and \ref{thm:rank stabilization} in Sections
\ref{sec:construction}, \ref{sec:lifts of maps} , and \ref{sec:pf of stable
thm} respectively.

\subsection{Bott Periodicity}\label{subsec:Bott per.}
In this subsection we show how Theorems \ref{thm:rank
stabilization}, \ref{thm:lifts of maps to configurations}, and Taubes'
stabilization
leads to an alternative, relatively quick proof of the following homotopy
equivalences in the periodicity spectrum:
\begin{thm}[Bott]\label{thm:Bott periodicity}
Let $SU$, $U$, $SO$, and $Sp$ denote the direct limit groups  of $SU(n)$,
$U(n)$, $SO(n)$, and $Sp(n)$ as $n\to \infty $. Let $\Omega ^{j}X$ denote
the $j$-fold loop space of $X$. The following are homotopy equivalences:
\begin{eqnarray*}
\Omega ^{2}SU&\sim &U,\\
\Omega ^{2}Sp&	\sim &	U/O,\\
\Omega ^{2}SO&	\sim &	U/Sp,\\
\Omega ^{4}SO&	\sim &	Sp,\\
\Omega ^{4}Sp&	\sim &	O.
\end{eqnarray*}
\end{thm}
\begin{rem}
The first equivalence is Bott periodicity for the unitary group and the
next four appear in the real periodicity spectrum. The only missing
homotopy equivalences:
$$
\Omega ^{2}(Sp/U)\sim BO\times \znums \text{  and  } \Omega ^{2}(SO/U)\sim
BSp\times \znums
$$
are related to monopoles (see Cohen and Jones \cite{Co-Jo} and the thesis
of Ernesto Lupercio \cite{Lupercio}).
\end{rem}

\proof
Let $i'$, $j'$ and $t'$ denote the maps on the moduli spaces induced by rank
inclusion, pull-back from $S^{4} $ to $\cpbar $, and Taubes' gluing
respectively ($i'$ and $j' $ are the descent of the maps $i$ and $j$ in
Theorem \ref{thm:lifts of maps to configurations}). We will argue that
$i'$, $j'$, and $t'$ commute up to homotopy. The maps  $i'$ and $j'$
commute (on the nose) from Theorem
\ref{thm:lifts of maps to configurations}; and we can see that $t'$
commutes up to homotopy with $i' $ and $j'$ from some general properties of
$t'$: The Taubes' map for any semi-simple compact Lie group $G$ is obtained
from the Taubes map for $SU (2)$ via any homomorphism $SU (2)\to G$
generating $\pi _{3} (G)$. Since the inclusions $G_{n}\embed G_{n'}$,
$n'>n$ ($n>4$ if $G_{n}=SO (n)$) induce isomorphisms on $\pi _{3}$, $t'$
automatically commutes with $i'$. To see that $t'$ commutes up to homotopy
with $j' $ we use almost instantons: connections with Yang-Mills energy
smaller than a small constant $\epsilon $. The space of almost instantons
$\M{k,\epsilon }{G_{n}} (X)$
has a strong deformation retract onto the space of instantons and there is
a map $t'_{\epsilon }:\M{k,\epsilon }{G_{n}} (X)\to \M{k+1,\epsilon }{G_{n}}
(X)$ homotopic to $t'$. It is local in the sense that $t_{\epsilon }'
(A)$ agrees with $A$ up to gauge in the complement of a ball about the
gluing point. On the other hand, the natural map $\cpbar \to S^{4}$ is a
conformal isometry on the complement of the hyperplane that gets mapped to a
point. Connections pulled back by this map have the same Yang-Mills
energy and we get a map $j'_{\epsilon }$ on almost instantons. Thus as long
as we choose our gluing point away from the hyperplane, $t'_{\epsilon }$
and $j'_{\epsilon }$ commute and so $t'$ and $j'$ commute up to homotopy.

The maps $t'$, $j'$, and $j$ then
induce commuting maps on the corresponding direct limit moduli and
configurations spaces when
$n\to \infty $ . We will assume that we have passed to that limit
throughout the rest of this section. From Theorem \ref{thm:rank
stabilization}, we have that
$A^{G_{\infty }}_{k}(X)$ is contractible. We can thus identify the homotopy
fibers
of the $j'$ maps to get the following fibrations:
\footnote{As we will see in Table 1, the structure groups of the various
principle bundles
$A^{G_{n}}_{k}(X)$ are the complex forms of the groups $U(k)$, $Sp(k/2)$, and
$O(k)$ (\ie $Gl(k,\cnums )$, $Sp(k/2,\cnums )$, and $O(k,\cnums )$). Since
the complex forms of the groups are homotopic to their compact form, their
classifying spaces are the same (up to homotopy). It
is traditional to use the compact form when refering to classifying spaces
so we will use the notation of the compact group for the rest of this
section. }
\begin{equation}\label{eqn:j-map fibrations}
\begin{CD}
U(k)\times U(k)/U(k)@>>>\M{k}{SU}(S^{4})@>{j'}>>\M{k}{SU}(\cpbar ),\\
U(k)/O(k)@>>>\M{k}{Sp}(S^{4})@>{j'}>>\M{k}{Sp}(\cpbar ),\\
U(k)/Sp(k/2)@>>>\M{k}{SO}(S^{4})@>{j'}>>\M{k}{SO}(\cpbar )
\end{CD}
\end{equation}
where $U(k)$ is included into $U(k)\times U(k)$ via the diagonal. Here we
are using the fact that Theorem \ref{thm:lifts of maps to configurations}
gives us $j$, the lift of $j'$ to the principle bundles $A^{G_{n}}_{k}
(X)$ that intertwines the actions.

Since $j'$ commutes with $t$, the above fibrations are valid for the direct
limit spaces when $k\to \infty $.

We now use Taubes' theorem to compare the above fibrations with the
fibration on the space of connections induced by the cofibration
$S^{2}\embed \cpbar \to S^{4}$. Let $\mathcal{B}^{G_{n}}_{k}(X)$ denote the
space of all $G_{n}$-connections of charge $k$ modulo based gauge
equivalence.  $\mathcal{B}^{G_{n}}_{k}(X)$ is homotopy equivalent to
the mapping space $Map_{k}(X,BG_{n})$ and the cofibration $S^{2}\embed
\cpbar \to S^{4}$ gives rise to a fibration
$$
\Omega _{k}^{4}BG_{n}\xrightarrow{j'} Map_{k}(\cpbar ,BG_{n})\to \Omega
^{2}BG_{n} .
$$
Up to homotopy, the map $j'$ in the above sequence is induced by pulling
back connections via $\cpbar \to S^{4}$ (thus justifying the notation).
These maps also commute with the group inclusions $i$ and so give a
fibration in the $n\to \infty $ limit. Also, $\mathcal{B}^{G_{n}}_{k}(X)$
and $\mathcal{B}^{G_{n}}_{k+1}(X)$ are naturally homotopy equivalent and so
we implicitly identify them and drop the notational dependence; we have:
\begin{equation}\label{eqn:mapping space fibration}
\mathcal{B}^{G}(S^{4})\xrightarrow{j'}\mathcal{B}^{G}(\cpbar )\to \Omega G.
\end{equation}

Taubes' stabilization theorem states that the inclusions
$\M{k}{G_{n}}(X)\embed \mathcal{B}^{G_{n}}(X)$ induce a homotopy
equivalence in the limit $k\to \infty $. Since the inclusion of the moduli
spaces into $\mathcal{B}$ commutes with both $j'$ and $i'$ we can pass to
the $k\to \infty $ and $n\to \infty $ limits and use Equation
\ref{eqn:j-map fibrations} to identify the homotopy fiber of
$j':\mathcal{B}^{G}(S^{4})\to\mathcal{B}^{G}(\cpbar ) $. This fiber is in
turn homotopy equivalent to $\Omega ^{2}G$ by the sequence \ref{eqn:mapping
space fibration}. Thus for $G=SU$, $Sp$, and $SO$ respectively, we get
\begin{eqnarray*}
\Omega ^{2}SU&	\sim &	U\times U/U\sim U,\\
\Omega ^{2}Sp&	\sim &	U/O,\text{ and}\\
\Omega ^{2}SO&	\sim &	U/Sp.
\end{eqnarray*}

The final two homotopy equivalences are arrived at by applying Theorem
\ref{thm:rank stabilization}  and Taubes' stabilization directly to
$\M{k}{Sp}(S^{4})$ and $\M{k}{SO}(S^{4})$ .
\qed

\subsection{The algebro-geometric moduli spaces and $\U
_{k}^{G_{n}}$}\label{subsec: alg-geo and U's}

>From now on, we will use the notation $X$ and $Y$ to denote $S^{4}$
and $\cnums \P ^{2} $
or $\cpbar$ and $\til{\cnums \P}^{2}$ (Recall from the introduction that
$\til{\cnums\P}^{2}$ is the blown-up projective plane). Consider the moduli
space $\M{alg}{n,k}(Y)$ consisting of pairs $(\E ,\tau )$ where $\E $ is a
rank $n$ holomorphic bundle on $Y$ and $\tau
:\E|_{H}\stackrel{\cong}{\longrightarrow} \cnums ^{n}\otimes \O _{H}$ is an
isomorphism.

Let $p:Y\to X $ be the smooth map that sends
$H\mapsto x_{0}$ and is one-to-one elsewhere. The map $p$ is compatible
with the natural orientations (we think of $p$ as a ``anti-holomorphic
blowdown''). We can construct a natural map
$$
\Xi :\M{k}{SU(n)}(X )\to \M{alg}{n,k}(Y)
$$
by defining the holomorphic structure on $p^{*}(E)$ corresponding to $\Xi
([A])$ to be $(d_{p^{*}A})^{0,1}$ and $\tau $ is induced by the fixed
isomorphism of $E_{x_{0}}$ (c.f. \cite{Bryan}).
\begin{thm}[Donaldson \cite{DonMonads}, King \cite{King}]\label{thm:King's thm}
The map $\Xi $ is an isomorphism of moduli spaces.
\end{thm}

Consider the moduli space $\M{alg,\pm }{n,k}(Y)$ of triples $(\E ,\tau
,\phi )$ where $(\E ,\tau )\in \M{alg}{n,k}(Y)$ and $\phi :\E \to \E ^{*}$
is a isomorphism such that $\phi ^{*}=\pm \phi $. We can construct maps
$\Xi _{\pm }$ in the same fashion as $\Xi $. A consequence of Theorem
\ref{thm:King's thm} is
\begin{cor}\label{cor:Sp and SO bundles from SU }
The maps
\begin{eqnarray*}
\Xi _{+}&:&\M{k}{SO(n)}(X )\to \M{alg,+}{n,k}(Y)\text{ and}\\
\Xi _{-}&:&\M{k}{Sp(n)}(X )\to \M{alg,-}{n,k}(Y)
\end{eqnarray*}
are moduli space isomorphisms ($(X,Y)$ is $(S^{4},\cnums \P ^{2})$
or $(\cpbar ,\til{\cnums \P }^{2})$).
\end{cor}
\proof

Recall that we consider $SO(n)$ or $Sp(n/2)$ connections to be $SU(n)$
connections that are compatible with a real or symplectic structure $\phi
$, \ie
$$
\nabla_{A^{*}}(\phi s)=\phi \nabla_{A}s
$$
where $A\in \mathcal{A(E)}$ and $A^{*}$ is the induced connection in
$\mathcal{A}(E^{*})$. Compatibility implies that $\phi $ will be a
holomorphic map with respect to the holomorphic structures defined by
$(d_{\pi^{*}A})^{0,1}$ and $(d_{\pi ^{*}(A^{*})})^{0,1}$. Conversely, let
$(\E ,\tau ,\phi )$ be in $\M{k}{Sp(n/2)}(X)$ or $\M{k}{SO(n)}(X)$. Choose a
hermitian structure on $\E $ compatible with $\phi $ and $\tau $. By
Theorem \ref{thm:King's thm}, the unique hermitian connection is the
pullback of an anti-self-dual $SU(n)$ connection on $E$ which is, by
construction, compatible with $\phi $. \qed

Henceforth we will drop the $\M{alg}{}$ notation and use $\M
{k}{G_{n}}(X )$ to refer to either moduli space.

The moduli space $\M{k}{G_{n}}(X) $ has a universal bundle (see
Lemma 3.2 of \cite{Br-Sa})
$$
\begin{CD}
\Bbb{E}\\
@VVV\\
\M{k}{G_{n}}(X)\times Y
\end{CD}
$$
so that $\Bbb{E}|_{\{\E  \}\times Y}\cong \E $.

Consider the cohomology groups $H^{i}(\E (-H))$. The fact that $\E $ is
trivial on $H$ (and thus on nearby lines) implies that $H^{i}(\E (-H))=0$
for $i=0$ or $2$ (see \cite{Br-Sa}). The Riemann-Roch theorem then gives
$\dim H^{1}(\E (-H))=k$. We will see from the construction of section
\ref{sec:construction} that the vector space $W$ of Table 1 can
be canonically identified with $H^{1}(\E (-H))$. In the case of $X=\cpbar
$ and $G_{n}=SU(n)$, a similar argument shows $\dim H^{1}(\E (-2H+E))=k$ and
$U$ is canonically identified with $H^{1}(\E (-2H+E))$.

Let $\pi :\M{k}{G_{n}}(X)\times Y\to \M{k}{G_{n}}(X)$ be projection. One
consequence of the above discussion is that the higher direct image sheaf
$R^{1}\pi _{*}(\Bbb{E}(-H))$ is a rank $k$ bundle on $\M{k}{G_{n}}(X)$.
Consequently, we have the following geometric interpretation of the
configuration spaces $A^{G_{n}}_{k}(X)$ (c.f. Appendix \ref{subsec:diff-geo
constr of univ bundles}):

\begin{thm}\label{thm:config space is principal bundle assoc to R1pi*}
The space of configurations $A^{G_{n}}_{k}(X) $ (see Table 1) is
homeomorphic to the total space of the frame bundle of $R^{1}\pi
_{*}(\Bbb{E}(-H))$ except for the case $A^{SU(n)}_{k}(\cpbar )$ which is
homeomorphic to the frame bundle of $$R^{1}\pi _{*}(\Bbb{E}(-H))\oplus
R^{1}\pi _{*}(\Bbb{E}(-2H+E)).$$
\end{thm}

\proof The fiber of $A^{G_{n}}_{k}(X)\to \M{k}{G_{n}}(X)$  over a point $\E
$ is the orbit of a
representative configuration by the automorphism group. This can be
identified with the space $\Iso (W,\cnums ^{k})$ (or $\Iso (W,\cnums
^{k})\times \Iso (U,\cnums ^{k})$ in the $X=\cpbar $, $G_{n}=SU(n)$ case)
where we also understand $\Iso (W,\cnums ^{k})$ to be isomorphisms of
symplectic or real vector spaces in the $X=S^{4}$, $G_{n}=SO(n)$ or
$Sp(n/2)$ cases.

\section{Construction of the moduli spaces} \label{sec:construction}
\subsection{Preliminaries}
To fill in Table 1, we rely heavily on the constructions of
Donaldson and King; we will recall what we need from their constructions in
subsection \ref{subsec:SU constructions}.
Let us first begin by introducing some general notation. For an
$n$-dimensional
projective manifold $M$ and a coherent sheaf $\mathcal{E}$ on $M$ let
$SD_{p,\mathcal{E}}$ denote the Serre duality isomorphism
$$
SD_{p,\mathcal{E}}:H^{p}(\mathcal{E})\to H^{n-p}(\E ^{*}(K))^{*}.
$$
Let $H^{i}(\phi ):H^{i}(\E\otimes \mathcal{G} )\to H^{i}(\mathcal{F}\otimes
\mathcal{G})$ denote the map in
cohomology induced by a sheaf map $\phi :\E \to \mathcal{F}$. If $s\in
H^{0}(\O (D))$ is a section vanishing on $D$, we let $\delta _{s}:H^{0}(\E
|_{D}) \to H^{1}(\E (-D))$ denote the coboundary map arising in the long
exact sequence associated to
$$
0\to \E (-D)\xrightarrow{s}\E \xrightarrow{r}\E _{D}\to 0.
$$
We will use the following elementary properties of Serre duality:

\begin{enumerate}
\item $SD_{p,\E }=(-1)^{p(n-p)}(SD_{n-p,\E ^{*}(K)})^{*}$,
\item $SD_{p,\E }$ is natural in the sense that
$$
\begin{CD}
H^{p}(\E )  @>{SD_{p,\E }}>> H^{n-p}(\E ^{*}(K))^{*}\\
@VV{H^{p}(\phi )}V   @VV{H^{n-p}(\phi ^{*})^{*}}V\\
H^{p}(\mathcal{F})@>{SD_{p,\mathcal{F}}}>> H^{n-p}(\mathcal{F^{*}}(K))^{*}
\end{CD}
$$
commutes.
\end{enumerate}

The sign in the first property arises from commuting the cup product.

A {\em monad} is a three term complex of
vector bundles over a complex manifold
$$\mathcal{U}
\xrightarrow{A}\mathcal{V}\xrightarrow{B}\mathcal{W}$$
such that $A$ is injective, $B$ is surjective, and
$B\comp A$ is $0$. The monad determines its {\em cohomology bundle}
$\mathcal{E}=\Ker(B)/\im(A)$.

The point is that one can build complicated
holomorphic bundles from relatively simple bundles using monads. By fixing
the bundles $(\mathcal{U} ,\mathcal{V},\mathcal{W})$ and allowing the maps
$A$ and $B$ to vary, one parameterizes a family of bundles. We say that
$(\mathcal{U},\mathcal{V},\mathcal{W})$ {\em effectively parameterizes
bundles} if the morphisms of $(\mathcal{U},\mathcal{V},\mathcal{W})$-monads
are in one-to-one correspondence with morphisms of the associated
cohomology bundles. This will be the case under favorable cohomological
conditions on $(\mathcal{U},\mathcal{V},\mathcal{W})$ (for details see
\cite{King}, \cite{Horrocks}).

The Chern character of the cohomology bundle can be computed by the formula
$$
ch(\E )=ch(\mathcal{V})-ch(\mathcal{U})-ch(\mathcal{W}).
$$

Note that the cohomology bundle associated to the dual monad
$$
\mathcal{U}^{*}\xrightarrow{B^{*}}\mathcal{V}^{*}
\xrightarrow{A^{*}}\mathcal{W}^*
$$
is the dual bundle $E^{*}$. We call a monad {\em self-dual } (or {\em
anti-self-dual} ) if it is of the form
$$
\mathcal{U}\xrightarrow{A}\mathcal{V}\xrightarrow{A^{*}\beta
^{*}}\mathcal{U}^{*}
$$
where $\beta :\mathcal{V}\to \mathcal{V}^{*}$ is a real (or symplectic)
structure; \ie $\beta ^{*}=\beta $ (or $\beta ^{*}=-\beta $).

A self-dual monad is isomorphic to its dual by the isomorphism
$(\1 ,\beta ,\1)$ and an anti-self-dual monad is isomorphic to its dual by
$(\1 ,\beta ,-\1 )$. Thus if
$(\mathcal{U},\mathcal{V},\mathcal{U}^{*})$ effectively parameterizes
bundles, then $(\1,\beta ,\pm \1 )$ induces a
real (or symplectic) structure $\phi :\E \to \E ^{*}$ on the cohomology
bundle.

\subsection{The $SU(n)$ constructions.} \label{subsec:SU constructions}
We wish to show that the $SU(n)$ configurations of Table 1 give
rise to bundles in $\M{k}{SU(n)}(X)$.

For $(a_{1},a_{2},b,c)\in A^{SU(n)}_{k}(S^{4})$ consider the sequence of
bundles on $Y=\cnums \P ^{2}$:
\begin{equation}\label{eqn:S4 monad sequence}
W\otimes \O (-H)\xrightarrow{A}(W\oplus W\oplus \cnums ^{n})\otimes \O
\xrightarrow{B}W\otimes \O (H)
\end{equation}
where
\begin{eqnarray*}
A &=&\left(\begin{array}{c}
x_{1}-a_{1}x_{3}\\
x_{2}-a_{2}x_{3}\\
cx_{3}
\end{array} \right),\\
B&=&\left(\begin{array}{ccc}
-x_{2}+a_{2}x_{3}&x_{1}-a_{1}x_{3}&bx_{3}
\end{array} \right)
\end{eqnarray*}
and $\langle x_{1 },x_{2},x_{3} \rangle$ generates $H^{1}(\O (H))$ and $H$
is the zero set of $x_{3}$.

The integrability condition is equivalent to $B\comp A=0$. We define
$A^{SU(n)}_{k}(S^{4})$ to be the open dense set of the integrable
configurations such that $A$ and $B$ are pointwise injective and
surjective respectively.
Thus for configurations in $A^{SU(n)}_{k}(S^{4})$, Sequence \ref{eqn:S4
monad sequence} is a monad. By computing Chern classes and restricting
Sequence \ref{eqn:S4 monad sequence} to $H$, one can see that the
cohomology bundle $\E $ lies in $\M{k}{SU(n)}(S^{4})$. In fact the converse
is true:
\begin{thm}[Donaldson]
Every $\E \in \M{k}{SU(n)}(S^{4})$ is given by a monad of the form in
Sequence \ref{eqn:S4 monad sequence} and the correspondence is unique up to
the natural action of $Gl(W)$. Furthermore, $W$ is canonically identified
with $H^{1}(\E (-H))$ (Okonek, et. al. \cite{oss} pg. 275).
\end{thm}

To finish the proof of Theorem \ref{thm:moduli construction} for $X=S^{4}$
and $G_{n}=SU(n)$ we only need to show that the automorphism group acts
freely on $A^{SU(n)}_{k}(S^{4})$. This follows from the identification of
$A^{SU(n)}_{k}(S^{4})$ with the frame bundle of $R^{1}\pi
_{*}(\Bbb{E}(-H))$ (see Theorem \ref{thm:config space is principal bundle
assoc to R1pi*}).

For  $(a_{1},a_{2},x,b,c)\in A^{SU(n)}_{k}(\cpbar )$ consider the sequence of
bundles on $Y=\til{\cnums \P} ^{2}$:
\begin{equation}\label{eqn:cpbar monad sequence}
\begin{array}{c}
{U\otimes \O (-H)}\\
\oplus\\
{W\otimes \O (-H+E)}
\end{array}
\stackrel{A}{\longrightarrow}
V\otimes \mathcal{O}\stackrel{B}{\longrightarrow}
\begin{array}{c}
{W\otimes \O (H)}\\
\oplus\\
{U\otimes \O (H-E)}
\end{array}
\end{equation}
where $V=W\oplus U\oplus W\oplus U\oplus \cnums ^{n}$ and
\begin{eqnarray*}
A&=&\left(\begin{array}{cc}
		a_1x_3 &-y_2\\
		x_1-xa_1x_3 & 0 \\
		a_2x_3 & y_1 \\
		x_2-xa_2x_3 & 0 \\
		cx_3 & 0
	  \end{array}\right)\label{eq:form of A}\\
B&=&\left(\begin{array}{ccccc}
		x_2 & a_2x_3 & -x_1 & -a_1x_3 & bx_3 \\
		xy_1 & y_1 & xy_2 & y_2 & 0
		\end{array}\right).\label{eq:form of B}
\end{eqnarray*}
We have chosen sections $\langle x_1,x_2,x_3 \rangle$ spanning
$H^0(\O (H))$
and $\langle y_1,y_2\rangle$ spanning $H^0(\O (H-E))$ so that $x_{3} $
vanishes on $H$ and $x_{1}y_{1}+x_{2}y_{2}$ spans the kernel of $H^{0}(\O
(H))\otimes H^{0}(\O (H-E))\to H^{0}(2H-E)$.

The integrability condition is equivalent to $B\comp A=0$. We define
$A^{SU(n)}_{k}(\cpbar )$ to be the open dense set of the integrable
configurations that are such that $A$ and $B$ are pointwise injective and
surjective respectively.
Thus for configurations in $A^{SU(n)}_{k}(\cpbar )$, Sequence
\ref{eqn:cpbar monad sequence} is a monad. By computing Chern classes and
restricting
Sequence \ref{eqn:cpbar monad sequence} to $H$, one can see that the
cohomology bundle $\E $ lies in $\M{k}{SU(n)}(\cpbar )$. Once again the
converse is true:
\begin{thm}[King]
Every $\E \in \M{k}{SU(n)}(\cpbar )$ is given by a monad of the form in
Sequence \ref{eqn:cpbar monad sequence} and the correspondence is unique up
to the natural action of $Gl(W)\times Gl(U)$. Furthermore, $W$ and $U$ are
canonically identified with $H^{1}(\E (-H))$ and $H^{1}(\E (-H+E))$
respectively.
\footnote{Historically, the
correspondence between holomorphic bundles and instantons on $S^{4} $ or
$\cpbar $ was proved by constructing the bundle moduli spaces as in this
section and showing that the construction is equivalent to the ``twistor''
construction of instantons. Now there is a direct analytic proof of the
correspondence due to Buchdahl \cite{Bu93} that also applies to $\cpbar \#
\cdots \# \cpbar$.}
\end{thm}
Once again we see that automorphism group acts
freely on $A^{SU(n)}_{k}(\cpbar )$ from the identification of
$A^{SU(n)}_{k}(\cpbar )$ with the frame bundle of $$R^{1}\pi
_{*}(\Bbb{E}(-H))\oplus R^{1}\pi_{*}(\Bbb{E}(-H+E)).$$

\begin{rem}\label{rem:dual configs}
If $\E \in \M{k}{SU(n)}(X)$ then $\E ^{*}\in \M{k}{SU(n)}(X)$ and is given
by the cohomology of the dual monad. To find the ``dual configuration'', we
need to use a monad automorphism to put the dual monads into the form
determined by a configuration. One can then see that the correspondence $\E
\mapsto \E ^{*}$ is realized on the level of configurations by
$$
(a_{1},a_{2}, b,c)\mapsto (a_{1}^{*},a_{2}^{*},-c^{*},b^{*})
$$
in the $S^{4}$ case and
$$
(a_{1},a_{2},x, b,c)\mapsto (a_{1}^{*},a_{2}^{*},x^{*},-c^{*},b^{*})
$$
in the $\cpbar $ case.
\end{rem}
 We also will need some finer information about these constructions.
Namely, there is cohomological interpretations for the maps occurring in
the configurations.
For $(a_{1},a_{2},b,c)\in A^{SU(n)}_{k}(S^{4})$ the maps are given by the
following compositions:
\begin{eqnarray*}
a_{i}&:&H^{1}(\E (-H))\xrightarrow{H^{1}(x_{3})^{-1}}H^{1}(\E
(-2H))\xrightarrow{H^{1}(x_{i})} H^{1}(\E (-H)),
\label{eqn:coh interp of ai in S4 case}\\
b&:&H^{0}(\E |_{H})\xrightarrow{\delta _{x_{3}}}H^{1}(\E (-H)),\\
c^{*}&:&H^{0}(\E ^{*}|_{H})\xrightarrow{\delta _{x_{3}}}H^{1}(\E
^{*}(-H))\xrightarrow{H^{1}(x_{3})^{-1}}H^{1}(\E
^{*}(-2H))\xrightarrow{SD}H^{1}(\E (-H))^{*}.
\end{eqnarray*}

 With our definition of $\langle x_{1},x_{2},x_{3}
\rangle$ and $\langle y_{1},y_{2} \rangle$ we get a well defined section
$s=x_{2}/y_{1}=-x_{1}/y_{2}$ of $H^{0}(\O (E))$.
For $(a_{1},a_{2},x,b,c)\in A^{SU(n)}_{k}(\cpbar )$ the maps are given by the
following compositions:
\begin{eqnarray*}
a_{1}&:&H^{1}(\E (-H))\xrightarrow{H^{1}(x_{3})^{-1}}H^{1}(\E
(-2H))\xrightarrow{H^{1}(-y_{2})} H^{1}(\E (-H)),\\
a_{2}&:&H^{1}(\E (-H))\xrightarrow{H^{1}(x_{3})^{-1}}H^{1}(\E
(-2H))\xrightarrow{H^{1}(y_{1})} H^{1}(\E (-H)),\\
x&:&H^{1}(\E (-H+E))\xrightarrow{H^{1}(s)}H^{1}(\E (-H)),\\
b&:&H^{0}(\E |_{H})\xrightarrow{\delta _{x_{3}}}H^{1}(\E (-H)),\\
c^{*}&:&H^{0}(\E ^{*}|_{H})\xrightarrow{\delta _{x_{3}}}H^{1}(\E
^{*}(-H))\xrightarrow{H^{1}(x_{3})^{-1}}H^{1}(\E
^{*}(-2H))\xrightarrow{SD}H^{1}(\E (-H+E))^{*}.
\end{eqnarray*}

King gives a detailed discussion of this description. In the $S^{4} $ case,
one can also ferret these maps out of the Beilinson spectral sequence
derivation of the monads on $\cnums \P^{2}$ (\cite{oss} pg. 249-251,275)
using the triviality of $\E $ on $H$.

\begin{rem}\label{rem:jumping lines}
In general, if $z$ is the defining section of a divisor $D$ and $D$ is
geometrically a rational curve, then it is easy to see from the long exact
sequence that
$$
H^{1}(z):H^{1}(\E (-2D))\to H^{1}(\E (-D))
$$
is an isomorphism if and only if $\E |_{D}$ is trivial. Thus, in the
above cohomological interpretations of configurations, $H^{1}(x_{3})$ is
always an isomorphism. We see then, for example, that the map $x$ is
singular if and only if $\E $ has the exceptional curve $E$ as a ``jumping
line''. Likewise we can interpret $a_{1}$ and $a_{2}$: The complement of
$H$ in $Y$ is either a complex plane or a complex plane blown-up at the
origin. In either case lines through the origin are given by the zeros of
$\mu _{1}x_{1}+\mu _{2}x_{2}$. Thus $\E $ will have jumping lines at those
lines parameterized by $(\mu _{1},\mu _{2})$ for which $\mu
_{1}a_{1}+\mu _{2}a_{2}$ is singular. This circle of ideas has been
utilized heavily by Hurtubise, Milgram, {\it et. al.} who use jumping lines
to give a filtration of the moduli spaces (\cite{Hurtubise},
\cite{Hur-Mil}, \cite{BHMM}).
\end{rem}

\subsection{Construction of $\M{k}{Sp(n/2)}(X)$ and $\M
{k}{SO(n)}(X)$}\label{subs: construction of MSP and MSO}

We now use the constructions of the previous subsection to construct the
moduli spaces $\M{k}{Sp(n/2)}(X)$ and $\M{k}{SO(n)}(X)$.
We first show that given $Sp(n/2)$ or $SO(n)$ configurations from Table
1, one produces an appropriate self-dual or anti-self-dual
monad determining  an element of the corresponding moduli space. We then
show the converse, \ie given an element $(\E ,\tau ,\phi )$ of
$\M{k}{Sp(n/2)}(X)$ or $\M
{k}{SO(n)}(X)$ we can get an equivalence class of the corresponding
configurations from Table 1.

For each of the four cases with $G_{n}=Sp(n/2)$ or $SO(n)$, we will use
configurations to define a sequence
$$
\mathcal{U}\xrightarrow{A}\mathcal{V}\xrightarrow{A^{*}\beta
^{*}}\mathcal{U^{*}} .
$$
For integrable configurations (those in $\overline{A}^{G_{n}}_{k}(X)$), the
sequence will satisfy
$$
A^{*}\beta ^{*}A=0
$$
and for each of the cases we define
$A^{G_{n}}_{k}(X)\subset \overline{A}^{G_{n}}_{k}(X)$ to be the open dense
subset such that the corresponding map $A$ is pointwise injective. The sequence
will then be a monad.

For $(\alpha _{1},\alpha _{2},\gamma )\in
{A}^{Sp(n/2)}_{k}(S^{4})$ we define an anti-self-dual
monad by
\begin{equation}\label{eqn:S4 Sp-monad sequence}
W\otimes \O (-H)\xrightarrow{A}(W\oplus W\oplus \cnums ^{n})\otimes \O
\xrightarrow{A^{*}\beta ^{*}}W^{*}\otimes \O (H)
\end{equation}
where
\begin{eqnarray*}
A &=&\left(\begin{array}{c}
x_{1}-\alpha _{1}x_{3}\\
x_{2}-\alpha _{2}x_{3}\\
\gamma x_{3}
\end{array} \right),\\
\beta &=&\left(\begin{array}{ccc}
0&	\Phi &	0\\
-\Phi &	0&	0\\
0&	0&	J
\end{array} \right).
\end{eqnarray*}

For $(\alpha _{1},\alpha _{2},\gamma )\in
{A}^{SO(n)}_{k}(S^{4})$ we define an self-dual monad by
\begin{equation}\label{eqn:S4 SO-monad sequence}
W\otimes \O (-H)\xrightarrow{A}(W\oplus W\oplus \cnums ^{n})\otimes \O
\xrightarrow{A^{*}\beta ^{*}}W^{*}\otimes \O (H)
\end{equation}
where
\begin{eqnarray*}
A &=&\left(\begin{array}{c}
x_{1}-\alpha _{1}x_{3}\\
x_{2}-\alpha _{2}x_{3}\\
\gamma x_{3}
\end{array} \right),\\
\beta &=&\left(\begin{array}{ccc}
0&	\Phi &	0\\
-\Phi &	0&	0\\
0&	0&	\1
\end{array} \right).
\end{eqnarray*}

For $(\alpha _{1},\alpha _{2},\xi ,\gamma )\in
{A}^{Sp(n/2)}_{k}(\cpbar )$ we define an anti-self-dual
monad by
\begin{equation}\label{eqn:cpbar Sp monad}
\begin{array}{c}
{W^{*}\otimes \O (-H)}\\
\oplus\\
{W\otimes \O (-H+E)}
\end{array}
\stackrel{A}{\longrightarrow}
V'\otimes \mathcal{O}\xrightarrow{A^{*}\beta
^{*}}
\begin{array}{c}
{W\otimes \O (H)}\\
\oplus\\
{W^{*}\otimes \O (H-E)}
\end{array}
\end{equation}
where
\begin{eqnarray*}
A&=&\left(\begin{array}{cc}
		\alpha _1x_3 &-y_2\\
		x_1-\xi \alpha _1x_3 & 0 \\
		\alpha _2x_3 & y_1 \\
		x_2-\xi \alpha _2x_3 & 0 \\
		\gamma x_3 & 0
	  \end{array}\right),\\
\beta &=&\left(\begin{array}{ccccc}
0&	0&	\xi &	1&	0\\
0&	0&	1&	0&	0\\
-\xi &	-1&	0&	0&	0\\
-1&	0&	0&	0&	0\\
0&	0&	0&	0&	J
\end{array}\right).
\end{eqnarray*}
and $V'=W\oplus W^{*}\oplus W\oplus W^{*}\oplus \cnums ^{n}$.

For $(\alpha _{1},\alpha _{2},\xi ,\gamma )\in
{A}^{SO(n)}_{k}(\cpbar )$ we define a self-dual
monad by
\begin{equation}\label{eqn:cpbar SO monad}
\begin{array}{c}
{W^{*}\otimes \O (-H)}\\
\oplus\\
{W\otimes \O (-H+E)}
\end{array}
\stackrel{A}{\longrightarrow}
V'\otimes \mathcal{O}\xrightarrow{A^{*}\beta
^{*}}
\begin{array}{c}
{W\otimes \O (H)}\\
\oplus\\
{W^{*}\otimes \O (H-E)}
\end{array}
\end{equation}
where
\begin{eqnarray*}
A&=&\left(\begin{array}{cc}
		\alpha _1x_3 &-y_2\\
		x_1-\xi \alpha _1x_3 & 0 \\
		\alpha _2x_3 & y_1 \\
		x_2-\xi \alpha _2x_3 & 0 \\
		\gamma x_3 & 0
	  \end{array}\right),\\
\beta &=&\left(\begin{array}{ccccc}
0&	0&	\xi &	1&	0\\
0&	0&	-1&	0&	0\\
-\xi &	-1&	0&	0&	0\\
1&	0&	0&	0&	0\\
0&	0&	0&	0&	1
\end{array}\right).
\end{eqnarray*}
and $V'=W\oplus W^{*}\oplus W\oplus W^{*}\oplus \cnums ^{n}$.

By computing the Chern characters and restricting the monads to $H$ one can
see that the cohomology bundle of Monads \ref{eqn:S4 Sp-monad sequence} and
\ref{eqn:S4 SO-monad sequence} lie in $\M{k}{SU(n)}(S^{4})$ and the
cohomology bundle of Monads \ref{eqn:cpbar Sp monad} and \ref{eqn:cpbar SO
monad} lie in
$\M{k}{SU(n)}(\cpbar )$. Furthermore, since the Monads \ref{eqn:S4 Sp-monad
sequence}  and \ref{eqn:cpbar Sp monad} are anti-self-dual they induce a
symplectic structure $\phi:\E \to \E ^{*} $ on the corresponding cohomology
bundle which restricts to $J$ on $H$. The Monads \ref{eqn:S4 Sp-monad
sequence}  and \ref{eqn:cpbar Sp monad} thus define elements of
$\M{k}{Sp(n/2)}(S^{4})$ and
$\M{k}{Sp(n/2)}(\cpbar ) $ respectively. Similarly, the Monads \ref{eqn:S4
SO-monad sequence} and \ref{eqn:cpbar SO monad} define elements of
$\M{k}{SO(n)}(S^{4})$ and $\M{k}{SO(n)}(\cpbar )$. Finally, the
group of monad automorphisms that preserve the given form of the above monads
is induced by the natural action of the configuration automorphism groups
listed in Table 1.

Now suppose that $(\E ,\tau ,\phi )$  is an element of
$\M{k}{Sp(n/2)}(X)$ or $\M
{k}{SO(n)}(X)$. We wish to produce an equivalence class of the corresponding
configurations. Let $(a_{1},a_{2},b,c)$ or  $(a_{1},a_{2},x,b,c)$ be
a representative  configuration for $(\E ,\tau )\in \M{k}{SU(n)}(X)$.

We begin by defining the map $\Phi $ by the following composition of
isomorphisms:
\begin{equation}\label{eqn:defn of Phi}
\Phi =H^{1}(\phi )^{*}\comp SD_{\E (-2H)}\comp H^{1}(x_{3})^{-1}.
\end{equation}
For $X=\cpbar $, $\Phi$ is a map from $W$ to $U^{*}$ and for $X=S^{4}$,
$\Phi$ is a map from $W$ to $ W^{*}$.

\begin{prop}\label{prop: commuting relations for Phi}
When $X=S^{4}$ the map $\Phi $ satisfies the following relations:
\begin{eqnarray*}
\Phi ^{*}&=&\begin{cases}
\Phi &	\text{when $G_{n}=Sp(n/2)$,}\\
-\Phi &\text{when $G_{n}=SO(n)$,}
\end{cases}\\
\Phi a_{i}&=&a_{i}^{*}\Phi,\\
c^{*}&=&\begin{cases}
\Phi bJ& \text{if $G_{n}=Sp(n/2)$},\\
\Phi b&\text{if $G_{n}=SO(n)$ }.
\end{cases}
\end{eqnarray*}
When $X=\cpbar $ the map $\Phi $ satisfies the following relations:
\begin{eqnarray*}
\Phi a_{i}&=&\begin{cases}
-a_{i}^{*}\Phi ^{*}&	\text{if $G_{n}=SO(n)$,}\\
a_{i}^{*}\Phi ^{*}&	\text{if $G_{n}=Sp(n/2)$,}
\end{cases}\\
x^{*}\Phi &=&\begin{cases}
-\Phi ^{*}x&	\text{if $G_{n}=SO(n)$,}\\
\Phi ^{*}x&	\text{if $G_{n}=Sp(n/2)$,}
\end{cases}\\
c^{*}&=&\begin{cases}
\Phi bJ& \text{if $G_{n}=Sp(n/2)$},\\
\Phi b&\text{if $G_{n}=SO(n)$ }.
\end{cases}
\end{eqnarray*}
\end{prop}

\proof The proof is a straight forward application of the properties of Serre
duality to the cohomological interpretation of the configuration maps and
the definition of $\Phi $. For example, if $X=S^{4}$ we have a
commutative diagram:
$$
\begin{CD}
H^{0}(\E ^{*}|_{H})@>{\delta _{x_{3}}}>>H^{1}(\E
^{*}(-H))@<{H^{1}(x_{3})}<<H^{1}(\E ^{*}(-2H))@>{SD}>>H^{1}(\E
(-H))^{*}\\
@AA{H^{0}(\phi |_{H})}A@AA{H^{1}(\phi )}A@AA{H^{1}(\phi
)}A@AA{H^{1}(\phi^{*} )^{*}}A\\
H^{0}(\E|_{H})@>{\delta _{x_{3}}}>>H^{1}(\E
(-H))@<{H^{1}(x_{3})}<<H^{1}(\E (-2H))@>{SD}>>H^{1}(\E^{*}
(-H))^{*}
\end{CD}
$$
Now follow the diagram from the lower left corner to the
upper right using both directions along the perimeter.  Since $\phi |_{H}$
is  $\1 $ in the $SO$ case, $J$ in the $Sp$ case, and $J^{-1}=-J$, we see
that $c^{*}\comp \phi |_{H}=\pm \Phi \comp b$ where $\phi ^{*}=\pm \phi
$. The result for $c^{*}$ follows and a similar  diagram shows the $\cpbar
$ case.

If $X=S^{4}$ we wish to show that $\Phi ^{*}=\mp \Phi $ when $\phi ^{*}=\pm
\phi $. Noting that $H^{1}(x_{3})^{*}=H^{1}(x_{3})$ we have the following
commutative diagram:
$$
\begin{CD}
H^{1}(\E (-2H))@>{H^{1}(x_{3})}>>{H^{1}(\E (-H))}@>{H^{1}(\phi )}>>
H^{1}(\E ^{*}(-H))\\
@VV{SD}V && @VV{SD}V\\
H^{1}(\E ^{*}(-H))^{*}@>{H^{1}(\phi ^{*})^{*}}>>H^{1}(\E
(-H))^{*}@>{H^{1}(x_{3})}>>H^{1}(\E (-2H))^{*}
\end{CD}
$$
Following the diagram clockwise from the upper middle spot to the lower
middle spot gives the map $-\Phi ^{*}$ since $SD=-SD^{*}$ in this case.
Following the diagram counterclockwise yields $\pm \Phi $ when $\phi
^{*}=\pm \phi $ and so we have that $\Phi^{*} =\mp\Phi $.

We can prove the relation $\Phi a_{i}=\mp a_{i}^{*}\Phi^{*}$ in a similar
fashion. We write the relations algebraically and suppress the diagram:

\begin{eqnarray*}
\Phi a_{i}&=&H^{1}(\phi )^{*}\comp SD\comp H^{1}(x_{3})^{-1}\comp
H^{1}(z_{i})\comp  H^{1}(x_{3})^{-1}\\
&=&\pm H^{1}(\phi ^{*})^{*}\comp SD \comp  H^{1}(x_{3})^{-1}\comp
H^{1}(z_{i})\comp  H^{1}(x_{3})^{-1}\\
&=&\pm SD\comp H^{1}(\phi )\comp  H^{1}(x_{3})^{-1}\comp
H^{1}(z_{i})\comp  H^{1}(x_{3})^{-1}\\
&=&\mp  H^{1}(x_{3})^{-1}\comp
H^{1}(z_{i})\comp  H^{1}(x_{3})^{-1}\comp SD^{*}\comp H^{1}(\phi )\\
&=&\mp a_{i}^{*} \Phi ^{*}
\end{eqnarray*}
where $z_{i}=x_{i}$ in the $S^{4}$ case and $(z_{1},z_{2})=(-y_{2},y_{1})$
in the $\cpbar $ case.

Finally, we also have
\begin{eqnarray*}
x^{*}\Phi &=&H^{1}(s)^{*}\comp H^{1}(\phi )^{*}\comp SD\comp
H^{1}(x_{3})^{-1} \\
&=&\pm (H^{1}(\phi ^{*})\comp H^{1}(s))^{*}\comp SD\comp H^{1}(x_{3})^{-1}\\
&=&\pm H^{1}(x_{3})^{-1}\comp SD \comp H^{1}(\phi )\comp H^{1}(s)\\
&=&\mp H^{1}(x_{3})^{-1}\comp SD^{*} \comp H^{1}(\phi )\comp H^{1}(s)\\
&=&\mp \Phi ^{*}x.
\end{eqnarray*}
\qed

We are now in a position to define the inverse construction producing
configurations from $(\E ,\tau ,\phi )$. Define
a $Sp$ or $SO$ configuration $(\alpha _{1},\alpha _{2},\gamma )$ on $S^{4}$
by $\alpha _{i}=a_{i}$ and $\gamma =c$.  Define
a $Sp$ or $SO$ configuration $(\alpha _{1},\alpha _{2},\xi ,\gamma )$ on
$\cpbar $ by $\alpha _{i}=a_{i}(\Phi ^{-1})^{*}$, $\xi =\Phi ^{*}x $ and
$\gamma =c(\Phi ^{-1})^{*}$.
The proposition then implies that these are integrable configurations. This
correspondence intertwines the action of the automorphism group and is well
defined on the quotient. It is also the inverse  to the monad construction
and so completes the proof of  Theorem \ref{thm:moduli construction}.

\section{Lifting of maps to configurations}\label{sec:lifts of maps}

In this section we define the maps $i$ and $j$ and prove Theorem
\ref{thm:lifts of maps to configurations}. They will be maps on
configurations that descend to the maps on the moduli spaces. The map $i$
will descend to the map induced by the inclusion $G_{n}\embed G_{n'}$ for
$n'>n$ and $j$ will descend to the map induced by pulling back connections
via the map $\cpbar \to S^{4}$. The maps will intertwine the action of the
automorphism
groups, \ie $i$ will be equivariant (the automorphism groups are
independent of the rank), and $j$ will intertwine the action with
natural inclusions of the appropriate automorphism groups.

\subsection{The rank inclusion map $i$.}
We define $i$ on the various kinds of configurations by:

\begin{eqnarray*}
i&:&(a_{1},a_{2},b,c)\mapsto (a_{1},a_{2},b',c')\\
i&:&(a_{1},a_{2},x,b,c)\mapsto (a_{1},a_{2},x,b',c')\\
i&:&(\alpha _{1},\alpha _{2},\gamma )\mapsto (\alpha _{1},\alpha
_{2},\gamma ')\\
i&:&(\alpha _{1},\alpha _{2},\xi ,\gamma )\mapsto (\alpha _{1},\alpha
_{2},\xi ,\gamma ')
\end{eqnarray*}
where $c'=(\begin{array}{c}0\\c\end{array})$, $b'=\left(\begin{array}{cc}0&
b\end{array} \right)$, $\gamma '=(\begin{array}{c}0\\ \gamma \end{array})$
and $0$ is the appropriate zero map to (or from) $\cnums ^{(n'-n)}$. The map
is obviously equivariant with respect to the
automorphism groups and from the monad constructions it is easy to see that
$i$ descends to the map $\E \mapsto \E \oplus \O ^{(n'-n)}$. In terms of
connections, this is the map $A\mapsto A\oplus \Theta $ where $\Theta $ is
the trivial connection on the rank $n'-n$ bundle and this map is the
natural one induced by the inclusion $G_{n}\embed G_{n'}$.

\subsection{The pullback map $j$.}
We define the map $j$ as follows. For $G_{n}=Sp(n/2)$ or $SO(n)$ let
$$
j:(\alpha _{1},\alpha _{2},\gamma )\mapsto ( \alpha _{1}(\Phi
^{-1})^{*},\alpha _{2}(\Phi ^{-1})^{*}, \Phi ^{*},\gamma (\Phi ^{-1})^{*}).
$$
This map intertwines the actions of the automorphism groups and the natural
inclusions $SO(W)\embed Gl(W)$ or $Sp(W)\embed Gl(W)$ so it descends to a
map on the moduli spaces.

For $G_{n}=SU(n)$ we have automorphism groups $Gl(W)$ and $Gl(W)\times
Gl(U)$.  Choose an isomorphism $\chi :W\to U$ so that we can define an
inclusion $Gl(W)\embed Gl(W)\times Gl(U)$ by $g\mapsto (g,\chi g\chi
^{-1})$. Define $j$ to be
$$
j:(a_{1},a_{2},b,c)\mapsto (a_{1}\chi^{-1} ,a_{2}\chi^{-1} ,\chi
,b,c\chi ^{-1}).
$$
We see that $j$ then intertwines the action of the automorphism groups with
the inclusion $Gl(W)\embed Gl(W)\times Gl(U)$ induced by $\chi $. It is
clear that $j$ commutes with $i$.
\begin{lem}
The map $j$ induces the pull-back map on bundles.
\end{lem}
\proof
 We will
proceed by (1) pulling back the $\cnums \P ^{2}$ monad defined by an
$S^{4}$-configuration to a $\til{\cnums \P }^{2}$ monad via the blow-down
map $\til{\cnums \P }^{2}\to \cnums  \P ^{2}$ ({\it c.f.} subsection
\ref{subsec: alg-geo and U's}); (2) we use $\chi $ and a direct sum of the
monad with an exact sequence to get an equivalent monad of the form of
sequence \ref{eqn:cpbar monad sequence}; then (3) we will use a monad
automorphism
to arrive at the monad defined by the $\cpbar $-configuration $(a_{1}\chi
^{-1},a_{2}\chi ^{-1},\chi ,b,c\chi ^{-1})$. The $Sp$ and $SO$ cases are
similar and we leave them to the reader.

(1) Since in our notation $\langle x_{1},x_{2},x_{3} \rangle$ and $\O (\pm
H)$ on $\cnums \P ^{2}$ pull back to  $\langle x_{1},x_{2},x_{3} \rangle$
and $\O (\pm H)$ on $\til{\cnums \P }^{2}$, the pull-back of sequence
\ref{eqn:S4 monad sequence} does not change notationally. We apply the
monad isomorphism $(\chi ,\chi \oplus \chi \oplus \1 ,\1 )$ to it to get
\begin{equation}
\O (-H)\xrightarrow{A_{1}}(W\oplus W\oplus \cnums ^{n})\otimes \O
\xrightarrow{B_{1}}W\otimes \O (H)
\end{equation}
where
\begin{eqnarray*}
A_{1} &=&\left(\begin{array}{c}
x_{1}-\chi a_{1}\chi ^{-1}x_{3}\\
x_{2}-\chi a_{2}\chi ^{-1}x_{3}\\
c\chi^{-1} x_{3}
\end{array} \right),\\
B_{1}&=&\left(\begin{array}{ccc}
-\chi ^{-1}x_{2}+a_{2}\chi ^{-1}x_{3}&\chi ^{-1}x_{1}-a_{1}\chi
^{-1}x_{3}&bx_{3}
\end{array} \right).
\end{eqnarray*}

(2) Since $y_{1}$ and $y_{2}$ do not vanish simultaneously on $\til{\cnums
\P }^{2}$, the sequence
$$
W\otimes \O
(-H+E)\xrightarrow{\left(\begin{array}{c}-y_{2}\\y_{1}\end{array} \right)}
(W\oplus W)\otimes \O \xrightarrow{(\begin{array}{cc}\chi y_{1}&\chi
y_{2}\end{array})} U\otimes \O (H-E)
$$
is exact. We can thus direct sum this sequence to the previous monad to
obtain a monad with the same cohomology bundle. We get
\begin{equation}
\begin{array}{c}
{U\otimes \O (-H)}\\
\oplus\\
{W\otimes \O (-H+E)}
\end{array}
\stackrel{A_{2}}{\longrightarrow}
V\otimes \mathcal{O}\stackrel{B_{2}}{\longrightarrow}
\begin{array}{c}
{W\otimes \O (H)}\\
\oplus\\
{U\otimes \O (H-E)}
\end{array}
\end{equation}
where $V=W\oplus U\oplus W\oplus U\oplus \cnums ^{n}$ and
\begin{eqnarray*}
A_{2}&=&\left(\begin{array}{cc}
		0 &-y_2\\
		x_1-\chi a_1\chi ^{-1}x_3 & 0 \\
		0 & y_1 \\
		x_2-\chi a_2\chi ^{-1}x_3 & 0 \\
		c\chi ^{-1}x_3 & 0
	  \end{array}\right),\\
B_{2}&=&\left(\begin{array}{ccccc}
		0 & -\chi ^{-1}x_{2}+a_2\chi ^{-1}x_3 &0&\chi ^{-1}x_{1}
-a_1\chi ^{-1} x_3 & bx_3 \\
		\chi y_1 &0 & \chi y_2 &0 & 0
		\end{array}\right).
\end{eqnarray*}

(3) Finally we use an automorphism to put the monad into the form of the
sequence \ref{eqn:cpbar monad sequence} . Recall that
$s=-x_{1}/y_{2}=x_{2}/y_{1}$ is a well defined section in $H^{0}(\O (E))$. The
automorphism we use is
$(\eta _{1},\eta _{2},\eta _{3})$ where
$$
\eta_{1}=\left(\begin{array}{cc}1&0\\-\chi ^{-1}s&1\end{array} \right)
\text{,  } \eta_{3}=\left(\begin{array}{cc}1&\chi ^{-1}s\\0&1\end{array}
\right) ,
$$

$$
\eta _{2}=\left(\begin{array}{ccccc}
1&	-\chi ^{-1}&	0&	0&	0\\
0&	1&	0&	0&	0\\
0&	0&	1&	-\chi ^{-1}&	0\\
0&	0&	0&	1&	0\\
0&	0&	0&	0&	1
\end{array} \right)
$$
and matrix multiplication shows that $A=\eta _{2}A_{2}\eta _{1}^{-1}$ and
$B=\eta _{3}B_{2}\eta _{2}^{-1}$ are exactly the monad maps defined by the
$\cpbar $-configuration $(a_{1}\chi ^{-1},a_{2}\chi ^{-1},\chi ,b,c\chi
^{-1})$.
\qed

\section{Proof of the Stabilization theorem}\label{sec:pf of stable thm}

In section we prove Theorem \ref{thm:rank stabilization}.

We need to show that $\lim_{n\to \infty }A^{G_{n}}_{k}(X)$ is contractible.
Since the $A^{G_{n}}_{k}$'s are all algebraic spaces and the inclusion maps
are algebraic, they admit triangulations compatible with the
maps. Thus $A^{G_{\infty }}_{k}$ inherits the structure of a CW-complex and
so it suffices to show that the its homotopy groups are all zero. To this
end we show that the inclusion
$$
i:A^{G_{n}}_{k}(X)\to A^{G_{2k+n}}_{k}(X)
$$
is null homotopic. The basic point is that in $A^{G_{2k}}_{k} (X)$ there
are configurations whose only non-zero monad data consists of the maps to
or from $\cnums ^{2k}$, in other words the data $a_{i},\alpha _{i},x,$ or
$\xi $  are all zero. We will fix such a configuration in each case and
show that the image of $A^{G_{n}}_{k} (X)$ in $A^{G_{n+2k}}_{k} (X)$
homotopes to the image of the fixed configuration.

\begin{lem}\label{lem:existance of S1 invariant instantons}
There are configurations of the form $(0,\dots ,0,b_{0},c_{0})\in A^{SU
(2k)}_{k} (X)$ and $(0,\dots ,0,\gamma _{0})\in A^{Sp (2k)}_{k} (X)$ or $A^{SO
(2k)}_{k} (X)$.
\end{lem}
\proof
The integrability and non-degeneracy conditions for $SU $ configurations
reduce to $b_{0}c_{0}=0$ with $c_{0}$ injective and $b_{0}$
surjective. This can be easily accomplished by having $c_{0}$ map
isomorphically onto the first $k$ factors of $\cnums ^{2k} $ and $b_{0}$ an
isomorphism on the
remaining $k$ factors. For  the $Sp$ and $SO$ cases we need a map $\gamma
_{0}$ such that $\gamma_{0} ^{*}J\gamma_{0} =0$ or $\gamma_{0}
^{*}\gamma_{0} =0$ respectively, and so that $\gamma_{0} $ is
injective. This is also easily done;
for example, in the $SO$ case choose an isomorphism $Q:W\to \cnums ^{k}$
and let $\gamma _{0}= (Q,\sqrt{-1}Q)$ . We remark that
configurations of this form correspond exactly to instantons on $S^{4}$ or
$\cpbar $ that are invariant under the natural $S^{1}$ action.
\qed

Fix configurations as in the above lemma and define a homotopy
$H_{t}:A^{G_{n}}_{k}(X)\to A^{G_{2k+n}}_{k}(X)$ by the following:

For $X=S^{4}$ and $G_{n}=SU(n)$
$$
H_{t}(a_{1},a_{2},b,c)=((1-t)a_{1},(1-t)a_{2},(\begin{array}{cc}tb_{0}&(1-t)b
\end{array}),\left(\begin{array}{c}tc_{0}\\(1-t)c \end{array}
\right)) .
$$
For $X=\cpbar $ and $G_{n}=SU(n)$
$$
H_{t}(a_{1},a_{2},x,b,c)=((1-t)^{2/3}a_{1},(1-t)^{2/3}a_{2},(1-t)^{2/3}x,
(\begin{array}{cc}tb_{0}&(1-t)b
\end{array}),\left(\begin{array}{c}tc_{0}\\(1-t)c \end{array}
\right)) .
$$
For $X=S^{4} $ and $G_{n}=Sp(n/2)$ or $SO(n)$
$$
H_{t}(\alpha _{1},\alpha _{2},\gamma )=((1-t)\alpha
_{1},(1-t)\alpha _{2},\left(\begin{array}{c}t\gamma _{0}\\ (1-t)\gamma
\end{array} \right) ) .
$$
For $X=\cpbar $ and $G_{n}=Sp(n/2)$ or $SO(n)$
$$
H_{t}(\alpha _{1},\alpha _{2},\xi ,\gamma )=((1-t)^{2/3}\alpha
_{1},(1-t)^{2/3}\alpha _{2},(1-t)^{2/3}\xi ,\left(\begin{array}{c}t\gamma
_{0}\\ (1-t)\gamma
\end{array} \right))
$$

Configurations in the image of $H_{t}$ are integrable and non-degenerate so
$H_{t} $ is a well defined homotopy from the inclusion $i$ to a constant
map. We can thus conclude that in the $n\mapsto \infty $ limit
$A^{G_{n}}_{k} (X)$ is contractible and Theorem \ref{thm:rank
stabilization} follows. \qed

\appendix
\section{Differentio-geometric construction of universal bundles}
\label{subsec:diff-geo constr of univ bundles}

In this appendix we describe a differentio-geometric construction of the
universal bundles $R^{1}\pi _{*}(\Bbb{E}(-H))$ and $R^{1}\pi
_{*}(\Bbb{E}(-H+E))$ directly using anti-self-dual connections.

One motivation for the construction is the following. Let $\mathcal{FRED}$
denote the space of Fredholm operators on some Hilbert
space. $\mathcal{FRED}$ is a classifying space for $K$-theory in the sense
that
$$
K (X)\cong[X,\mathcal{FRED}]
$$
and the isomorphism is given by the index construction: if $f:X\to
\mathcal{FRED}$ is a family of Fredholm operators, then $[\Ker f
(x)]-[\Coker f (x)]$ pieces together to form a well defined element of $K
(X)$. On the other hand, $BU (k)$ classifies rank $k$ vector bundles:
$$
Vect_{k} (X)\cong [X,BU (k)]
$$
and one can ask if there is a natural, geometrically defined subset of
$\mathcal{FRED}$ that is homotopic
to $BU (k)$ such that the index construction induces the above
equivalence. By coupling instantons to the Dirac operator and using our
stabilization theorem we get such a family.

We define a rank $k$ bundle over $\M{k}{SU(n)}(S^{4})$ as follows: Since
$S^{4}$ is spin there are spin bundles $S^{\pm }$  and the Dirac operator
$\dirac $. We can couple the direct operator to a connection $A\in
\M{k}{SU(n)}(S^{4})$ to obtain an operator
$$
\dirac_A:\Gamma (S^{+}\otimes E)\to \Gamma (S^{-}\otimes E).
$$
Since the $\hat{A}(S^{4})=1$, the Atiyah-Singer index theorem shows that
the index of $\dirac _A$ is $-k$. The Bochner-Weitzenboch formula
for coupled Dirac operators is
$$
\dirac _{A}^{*}\dirac _{A}=\nabla_{A}^{*}\nabla_{A}+\frac{s}{4}+F_{A}^{+}.
$$
Since $F_{A}^{+}=0$ and $S^{4}$ has positive scalar curvature $s$, the
right hand side of the above equation is a positive operator. From the
usual argument we deduce that $\Ker (\dirac _{A})=0$ and so the vector
space $\mathfrak{W}_{A}=\Coker (\dirac _{A})$ is always $k$-dimensional.
The vector spaces $\mathfrak{W}_{A}$ vary smoothly with $A$ and piece
together to form the rank $k$ vector bundle $\mathfrak{W}\to
\M{k}{SU(n)}(S^{4})$ we
seek. The construction is natural with respect to the inclusion
$SU(n)\embed SU(n+1)$ and so $\mathfrak{W}$ ascends to $n\mapsto \infty $
direct limit. At the end of the section we will outline an argument showing
this is the same bundle as $R^{1}\pi _{*}(\Bbb{E}(-H))$.

For $X=\cpbar $, we consider the spin$_{\cnums }$ structure $W^{\pm }$
with $L=\det (W^{\pm })$ such that $c_{1} (L)$ is a generator of
$H^{2}(\cpbar ;\znums )$. Since $b^{+}_{2}(\cpbar )=0$ there is a
unique (up to gauge) connection $a\in \mathcal{A}(L)$
such that $F_{a}^{+}=0$. Using the connection $a$ we get a spin$_{\cnums }$
Dirac operator $\dirac _{a}$ which we can couple to a connection $A\in
\M{k}{SU(n)}(\cpbar )$ to get an operator
$$
\dirac _{a,A}:\Gamma (W^{+}\otimes E)\to \Gamma (W^{-}\otimes E).
$$

The Bochner-Weitzenboch formula for this operator is
$$
\dirac _{a,A}^{*}\dirac
_{a,A}=\nabla_{a,A}^{*}\nabla_{a,A}+\frac{s}{4}+F_{A}^{+}
+\frac{F_{a}^{+}}{2}.
$$
Once again the right hand side of this equation is positive so $\Ker
(\dirac _{a,A})=0$. We compute the index of the operator:
\begin{eqnarray*}
\operatorname{Ind}_{\cnums }(\dirac _{a,A})&=&ch(E)\cdot ch(L/2)\cdot
\hat{A}(\cpbar )[\cpbar ]\\
&=&(n-c_{2}(E))(1+c_{1}(L/2)+\tfrac{1}{2}c_{1}(L/2)^{2})
(1-\tfrac{1}{24}p_{1})[\cpbar]\\
&=&-k+n\left(\frac{c_{1}(L)^{2}-\sigma}{8} \right)\\
&=&-k.
\end{eqnarray*}
We define then a rank $k$ bundle $\mathfrak{W}_{+}\to \M{k}{SU(n)}(\cpbar
)$ whose fiber over $A$ is $\Coker (\dirac _{a,A})$. We get a second rank
$k$ bundle by the same construction with the spin$_{\cnums }$ structure
associated to $-c_{1}(L)$. We will not provide a complete proof that these
bundles are $R^{1}\pi _{*}(\Bbb{E}(-H))$ and  $R^{1}\pi
_{*}(\Bbb{E}(-H+E))$ but we will outline the argument. The case for $S^{4}$
is discussed in Section 3.3.3 and 3.3.4 in Donaldson and Kronheimer
\cite{D-K} so we will focus on the $\cpbar $ case.

Let $\til{\cnums}^{2}$ denote the complex plane blown-up at the origin with
the standard K\"ahler structure. $\til{\cnums }^{2}$ is biholomorphic to
$\til{\cnums \P }^{2}-H$ and conformally equivalent to $\cpbar -x_{0}$. We
wish to bring the operators $\dirac^{*} _{a,A}$ and $\dbar _{\E
}+\dbar ^{*}_{\E }$  to
$\til{\cnums }^{2}$ in order to compare them. Care must be taken on the
non-compact manifold to impose the correct decay conditions.

The conformally changed Dirac  operator is
$\til{\dirac}^{*}_{A,a}=e^{3f/2}\dirac^{*} _{A,a}
e^{-3f/2}$ where the conformal factor is $f=(1+r^{2}) $ and $r$ is the
radial coordinate in $\til{\cnums }^{2}$. Thus the kernel
of $\til{\dirac }^{*}_{A,a} $ consists of sections of order $\O (r^{-3})$
and this turns
out to be the same as $L^{2}$ harmonic spinors on $\til{\cnums }^{2}$. On
$\til{\cnums }^{2}$ the Dirac operator is $\dbar +\dbar ^{*}$ on $\Omega
^{01} ( E\otimes (K\otimes L)^{1/2})$. One can then identify the kernel
of $\dbar +\dbar ^{*}$ on $L^{2} $ section-valued $(0,1)$-forms on
$\til{\cnums }^{2}$ with the kernel of $\dbar +\dbar ^{*}$ on $\til{\cnums
\P }^{2}$ where the bundle has been twisted by a certain multiple of $H$ (the
``divisor at infinity'') determined by the decay conditions. In our case
the result is that $\Coker (\dirac^{*}_{A,a} )=\Ker (\dirac
_{A,a})=\Ker (\til{\dirac }_{A,a})=\Ker (\dbar +\dbar ^{*}) $ is $H^{1}(\E
(-2H+E))$ or  $H^{1}(\E(-2H))$  depending on the sign of $L$.


\end{document}